\documentclass[10pt,twocolumn,twoside]{IEEEtran}
\usepackage{cite}
\usepackage{amsmath,amssymb,amsfonts}
\usepackage{algorithmic}
\usepackage{graphicx}
\usepackage{textcomp}
\usepackage[switch]{lineno}

\usepackage{bm}
\usepackage[caption=false,font=normalsize,labelfont=sf,textfont=sf]{subfig}
\usepackage{textcomp,mathcomp}
\newtheorem{definition}{Definition}

\usepackage[some]{background}
\SetBgScale{1}
\SetBgContents{\parbox{20cm}{%
  \Huge Draft:  \today\\[18cm]\rotatebox{180}{\Huge Draft:  \today}}}
\SetBgColor{red}
\SetBgAngle{270}
\SetBgOpacity{0.2}
\usepackage[hidelinks]{hyperref}

\def\BibTeX{{\rm B\kern-.05em{\sc i\kern-.025em b}\kern-.08em
    T\kern-.1667em\lower.7ex\hbox{E}\kern-.125emX}}

\begin{document}



\title{CarboNet: A Finite-Time Combustion-Tolerant Compartmental Network for Tropospheric Carbon Control}

\author{Federico Zocco$^{1}$, Wassim M. Haddad$^{2}$, and Monica Malvezzi$^{1}$
\thanks{$^{1}$F. Zocco and M. Malvezzi are with the Department of Information Engineering and Mathematics, University of Siena, 53100 Siena, Italy
        {\tt\small federico.zocco.fz@gmail.com 
        monica.malvezzi@unisi.it}}%
\thanks{$^{2}$W. M. Haddad is with the School of Aerospace Engineering, Georgia Institute of Technology, Atlanta, GA 30332, USA.
        {\tt\small wassim.haddad@aerospace.gatech.edu}}%
\thanks{Corresponding author: Federico Zocco}
}

\maketitle

\begin{abstract}
While governments and international organizations have set the net-zero target to prevent a climate event horizon, practical solutions are lacking mainly because of the impracticability in completely replacing combustion processes. To address the net-zero target problem, in this paper we first design a compartmental network whose states must remain in the nonnegative orthant for physical consistency and in which the carbon dioxide emissions result from the combustion of diesel in vehicles and gas in house heaters. Then, we design both full-state and output-feedback linear-quadratic regulators of the compartmental network to bring the mass of carbon dioxide to the pre-industrial era, which is reached in approximately 25 and 60 days, respectively. The output-feedback controller tolerates for 6 days the combustion taking place in 5,000 vehicles and in 10,000 house heating systems, meets the net-zero target, and nullifies the extraction of finite natural resources. With closed-loop control, the tropospheric temperature stabilizes approximately to the pre-industrial era reference condition, i.e., to 13.5 °C, which is 21.7 °C lower than the steady-state temperature achieved without carbon capture. This work is a first step in designing optimal network control systems for climate stability. Source code is publicly available\footnote{\url{https://github.com/ciroresearch/CarboNet}}.         
\end{abstract}


\section{Introduction}\label{sec:introduction}
While governments and international organizations are adopting measures to prevent a climate event horizon, the temperature anomaly with respect to the pre-industrial era keeps increasing as recorded by NASA \cite{NASAtemperature}. The main cause of this phenomenon are the emissions of greenhouse gases, especially carbon dioxide \cite{NASAcarbon}. Carbon dioxide emissions result from combustion processes that are still the backbone of our economy despite the enormous international effort to adopt combustion-free energy and transportation systems, such as wind turbines and electric vehicles. Specifically, the share of primary energy consumption from fossil fuels in 2024 was approximately 73\% in Europe and the UK and 80\% in China and the United States \cite{DatafossilEnergy}, while the share of electric cars in use in 2024 was 4.7\% in Europe, 6.4\% in the UK, 11\% in China, and 2.7\% in the United States \cite{DataElecVehicles}.  

To the best of our knowledge, the analysis and control of Earth's climate via the methodological tools of dynamical systems and control theory is not widely explored as visible by searching for IEEE journal publications: as of 30 July 2025, the authors found only the works of Cavraro \cite{cavraro2024advocating}, Khargonekar \emph{et al.} \cite{khargonekar2024climate}, and De Nardi \emph{et al.} \cite{de2025climate}. A recent tutorial on control theory for climate science was recently carried out by Elsherif and Taha \cite{elsherif2025climate}, which highlights the lack of contributions from the control engineering community to the Earth's climate stabilization. Recently, Zocco \emph{et al.} \cite{zocco2025circularCarbon} began addressing carbon regulation as a network problem combining nonlinear control and reinforcement learning.

In this context, this paper makes the following contributions.
\begin{itemize}
\item{Developes a compartmental network model for carbon emissions and fuel consumption of a tunable number of vehicles and house heaters.}
\item{Designes full-state and output feedback linear-quadratic regulators (LQRs) for the network to bring the carbon dioxide level in a target urban area to the pre-industrial reference value. While carbon control has been already addressed in Zocco \emph{et al.} \cite{zocco2025circularCarbon} as a network control problem, this is the first paper where the control is centralized instead of distributed at different nodes.}
\item{Investigates the nonlinear dynamics of the tropospheric temperature in the target urban area for the closed-loop network.}
\item{Measures the circularity of the network to assess the consumption of natural resources and to insert the net-zero problem within the broader context of a circular economy.}
\end{itemize}

The following notation is adopted throughout the paper. Bold lower-case and capital letters are used for vectors and matrices, respectively; the set of nonnegative real numbers is denoted by $\overline{\mathbb{R}}_+$, a vector $\bm{x} \in \mathbb{R}^n$ with all nonnegative components is denoted by $\bm{x} \geq \geq 0$ or, equivalently, by $\bm{x} \in \overline{\mathbb{R}}^n_+$; similarly, a time-variant vector $\bm{x}(t)$ with all nonnegative components for all time is denoted by $\bm{x}(t) \in \overline{\mathbb{R}}^n_+$, $t \geq 0$.  

The remainder of the paper is organized as follows. Section \ref{sec:RelatedWork} discusses related work, while Section \ref{sec:LQRtheoryAndKeyVariables} covers the framework of linear-quadratic regulator designs along with the definition of the key variables we aim to regulate. Numerical examples and the case of $\text{CO}_2$ control are detailed in Section \ref{sec:NumStudies-fullStateFdbk} for full-state feedback control and Section \ref{sec:NumStudies-outputFdbk} for output feedback control. Finally, in Section \ref{sec:Conclusion} we provide conclusions.

\section{Related Work}\label{sec:RelatedWork}
\subsection{Circular Economy and Net Zero}
The net-zero target required by the United Nations \cite{UN-netZero} and the European Commission \cite{EU-netZero} requires to drastically reduce the carbon dioxide emissions. Net-zero essentially means that the amount of $\text{CO}_2$ entering the atmosphere is equal to the amount that is removed from it \cite{NZdefinition}. 

The circular economy paradigm is gaining attention as a solution to decouple economic growth from sustainable use of natural resources \cite{alexander2023handbook}. The practices of a circular economy include reducing the material demand, reusing material for as long as possible, and repairing products \cite{potting2017circular,velenturf2021principles}. These practices can reduce the carbon dioxide emissions by replacing the raw material extraction, transportation, and manufacturing with the less energy-intensive repair and reuse of functioning product parts and materials \cite{bonsu2020towards,yu2021evaluating,marche2022comparison}. Hence, increasing circularity helps meet the net-zero target. 

The \emph{circulation} of tropospheric carbon dioxide can further help meet the net-zero target mentioned above since carbon circulation yields a removal of carbon from the troposphere. Furthermore, carbon can be \emph{reused} through photosynthesis if plants are used for carbon sequestration as discussed in Section \ref{sub:CarbonCaptureSys} and in \cite{zhou2017bio,tsai2017potential,pereira2024nature}. In this paper, we increase the circularity of tropospheric carbon dioxide via properly controlled carbon-capture systems.

\subsection{Carbon Capture Systems}\label{sub:CarbonCaptureSys}
An intuitive approach to reduce the concentration of carbon dioxide in the troposphere is removing it. To this aim, different types of technologies are proposed in the literature. Rubin \emph{et al.} \cite{rubin2012outlook} analyzed three major solutions, namely, pre-combusion capture, post-combustion capture, and oxy-combustion capture. In power plants, the captured carbon dioxide can be piped off-shore or injected in underground locations with proper geological properties \cite{wilberforce2021progress}. The storage of carbon dioxide in fractured rocks was investigated by Romano \emph{et al.} \cite{romano2025evaluation}.

Nature-based solutions to carbon capture also exist and are divided into two macro-groups. The first group comprises of terrestrial carbon sinks such as forests, urban parks, green roofs, and green facades, while the second group comprises of water bodies such as lakes and rivers \cite{pereira2024nature}. In particular, the effectiveness of urban parks in the Italian city of Rome was deepened by Gratani \emph{et al.} \cite{gratani2016carbon}.

Microalgae are attracting attention because of their high efficiency in carbon absorption compared to other plants and because of their use as biomass at the end of their life \cite{zhou2017bio,tsai2017potential}. However, most of algal studies are currently limited to laboratory conditions \cite{yang2024carbon}. In this paper, we design a network involving a carbon sink and find an optimal absorption rate of carbon to recover the pre-industrial era reference condition. We do not assume any particular carbon capture system to perform the optimal absorption rate.

\subsection{Linear-Quadratic Regulators}
The linear-quadratic regulator (LQR) is a prominent approach for optimal control. It optimizes a positive quadratic cost function and it stabilizes linear time-invariant as well as time-varying systems to the origin \cite{nersesov2004optimal,ilka2022novel}. Bemporad \emph{et al.} \cite{bemporad2002explicit} addressed the case of constrained system inputs and outputs, Shi \emph{et al.} \cite{shi2012finite} developed finite-time horizon LQR with limited controller-system communication, while Maghfiroh \emph{et al.} \cite{maghfiroh2022improved} improved the LQR performance using particle swarm optimization and the Kalman filter. 

Network-control based on LQR was proposed by Jaleel and Shamma \cite{jaleel2017design}, in which each agent computes their control actions in real-time, by Duan \emph{et al.} \cite{duan2022distributed} where multi-input linear systems only need to communicate with their neighbors instead of with the whole network, and by Gao \emph{et al.} \cite{gao2021subspace} by leveraging the theory of graphons.

The major limitation of LQRs is that the choice of the weighting matrices is often based on empirical considerations and is problem dependant. To aid in this endeavor, approaches for the tuning of the weighting matrices were developed by Chacko \emph{et al.} \cite{chacko2024optimizing} and by Masti \emph{et al.} \cite{masti2021tuning}.      

LQRs can be designed not only in the case of access to the whole state, but also in the case of partial knowledge of the state. Solutions to the latter situation were provided by Levine and Athans \cite{levine2003determination}, by Nersesov \emph {et al.} \cite{nersesov2004optimal}, and by Ilka and Murgovski \cite{ilka2022novel}. 

Examples of applications of the LQR are wind turbine vibration control \cite{wang2024filter}, control of a single-phase inverter \cite{arab2019lqr}, and optimal drug delivery \cite{nersesov2004optimal}. In this paper, we develop both output and full-state feedback LQRs for network-based carbon control.

\section{Network, Key Variables, and Linear-Quadratic Regulators}\label{sec:LQRtheoryAndKeyVariables}
This section is divided into three parts. In the first part, the design of the compartmental network is covered, then the key sustainability variables to be controlled are introduced, and finally, the LQR design procedure is detailed. 
\subsection{Network Design}
\begin{figure}
\centering
\includegraphics[width=0.45\textwidth]{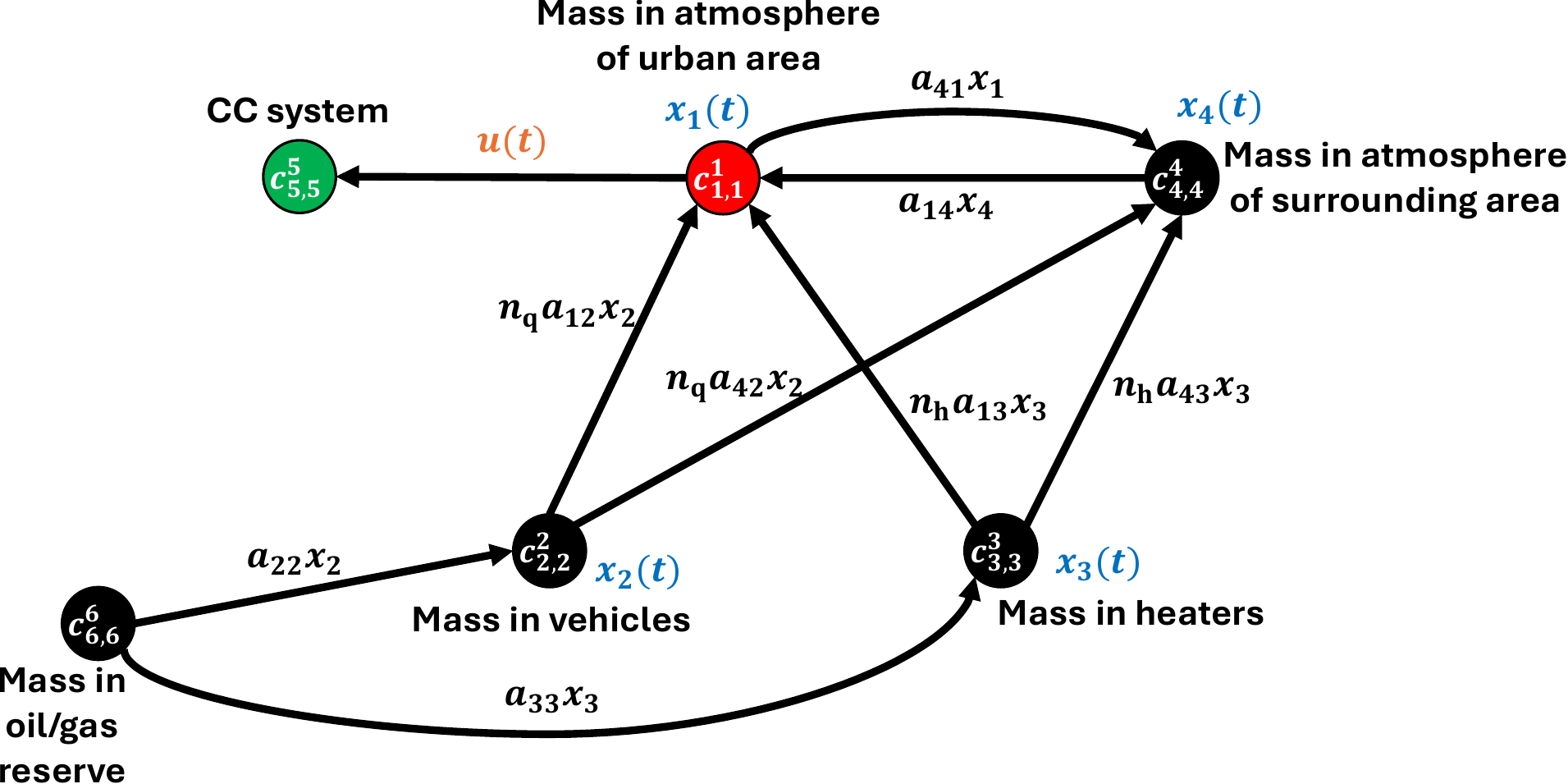}
\centering
\caption{Compartmental diagraph of the system for tropospheric carbon regulation. It can be represented as the set of compartments $\mathcal{N}_\text{c}$ (\ref{eq:setNc}).}
\label{fig:compDiagraph}
\end{figure}

The compartmental diagraph for carbon control considered in this paper is depicted in Fig. \ref{fig:compDiagraph} and it is based on the formalism of thermodynamical material networks \cite{zocco2023thermodynamical,zocco2025synchronized,zocco2025cirl,zocco2022circularity,zocco2025circular}. Specifically, the green node ($c^5_{5,5}$) is the carbon capture (CC) system. Its role is to remove the $\text{CO}_2$ from the red node ($c^1_{1,1}$), which indicates the urban area whose carbon dioxide we aim to regulate via the control laws covered in Section \ref{sec:LQRtheoryAndKeyVariables}. The tropospheric $\text{CO}_2$ in the urban area is increased by two types of sources: by $n_\text{q}$ vehicles with internal-combustion engines ($c^2_{2,2}$) and by $n_\text{h}$ heaters ($c^3_{3,3}$); the sources of carbon take and burn the fuel from the natural reserve of oil or gas ($c^6_{6,6}$). The urban area exchanges carbon with the surrounding area, which is indicated by $c^4_{4,4}$. Overall, the system has the four states $x_i(t)$ indicated in blue in Fig. \ref{fig:compDiagraph}, with $i \in \{1,2,3,4\}$, and the control signal $u(t)$ is indicated in orange. The states $x_1(t)$ and $x_4(t)$ are masses of $\text{CO}_2$ in the troposphere, while the states $x_2(t)$ and $x_3(t)$ are masses of fuel converted into $\text{CO}_2$ via combustion to move the $n_\text{q}$ vehicles (in $c^2_{2,2}$) or to heat-up the $n_\text{h}$ buildings (in $c^3_{3,3}$).        

The system in Fig. \ref{fig:compDiagraph} can be represented as the thermodynamical material network (see \cite{zocco2023thermodynamical,zocco2025synchronized,zocco2025cirl,zocco2022circularity,zocco2025circular})
\begin{equation}\label{eq:setNc}
\begin{aligned}
\mathcal{N}_c = \{c^1_{1,1},c^2_{2,2},c^3_{3,3},c^4_{4,4},c^5_{5,5},c^6_{6,6},c^7_{2,1},\\ c^8_{3,1},c^9_{4,1},c^{10}_{1,4},c^{11}_{1,5},c^{12}_{2,4},c^{13}_{3,4},c^{14}_{6,2},c^{15}_{6,3}\}.
\end{aligned}
\end{equation}
The arcs (i.e., the arrows) in the compartmental diagraph indicate the mass flow rates between the compartments, while their weights (i.e., their magnitudes) result from the product between a state $x_i(t)$ of the compartment $c^i_{i,i}$ and a rate constant $a_{k,s} \geq 0$, where for the pair ($k,s$) it holds that
\begin{equation}
(k,s) = 
\begin{cases}
\parbox[t]{.16\textwidth}{$k$ of the arc head $c^k_{k,k}$ and $s$ of the arc tail $c^s_{s,s}$,} & \quad\parbox[t]{.15\textwidth}{if the head and tail have states $x_k(t)$ and $x_s(t)$, respectively;} \\
\\
k = s, & \parbox[t]{.17\textwidth}{if at least one between the arc tail $k$ or the arc head $s$ is not associated with a state.} \\  
\end{cases}
\end{equation}
For the network (\ref{eq:setNc}), the number of compartments, the number of nodes, and the number of arcs are $n_\text{c} = 15$, $n_\text{v} = 6$, and $n_\text{a} = 9$, respectively. 

Thus, the state vector is $\bm{x}(t)$ = [$x_1(t)$, $x_2(t)$, $x_3(t)$, $x_4(t)$]$^\top \in \mathbb{\overline{R}}_+^4$, with the states indicated in Fig. \ref{fig:compDiagraph}, and the control input $u(t) \in \mathbb{R}$. Note that $\bm{x}(t) \geq \geq 0$, while $u(t) < 0$ indicates that the flow has a direction opposite to the one in Fig. \ref{fig:compDiagraph}. The application of the principle of mass balance to the continuous-time system in Fig. \ref{fig:compDiagraph} yields the following linear ordinary differential equations:
\begin{gather}
\begin{aligned}\label{eq:dynamicEq1}
\dot{x}_1(t) = -a_{41} x_1(t) + a_{14} x_4(t) + n_\text{h} a_{13} x_3(t) \\ + n_\text{q} a_{12} x_2(t) - u(t), \,\, x_1(0) = x_{10}, \,\, t \geq 0,
\end{aligned} 
\\
\begin{aligned}\label{eq:dynamicEq2}
\dot{x}_2(t) = a_{22}x_2(t) - n_\text{q} a_{12} x_2(t) - n_\text{q}a_{42}x_2(t), \\ x_2(0) = x_{20}, 
\end{aligned}
\\
\begin{aligned}\label{eq:dynamicEq3}
\dot{x}_3(t) = a_{33} x_3(t) - n_\text{h} a_{13} x_3(t) - n_\text{h} a_{43} x_3(t), \\ x_3(0) = x_{30},
\end{aligned}
\\
\begin{aligned}\label{eq:dynamicEq4}
\dot{x}_4(t) = a_{41} x_1(t) - a_{14} x_4(t) + n_\text{q} a_{42} x_2(t) \\ + n_\text{h} a_{43} x_3(t), \,\, x_4(0) = x_{40},
\end{aligned}
\end{gather}
which in state-space form is given by
\begin{gather}\label{eq:AxBu}
\dot{\bm{x}}(t) = \bm{A}\bm{x}(t) + \bm{B}u(t), \,\, \bm{x}(0) = \bm{x}_0, \,\, t \geq 0,
\\
\bm{y}(t) = \bm{C}\bm{x}(t),
\end{gather}
where $\bm{x}(t) \in \mathbb{\overline{R}}^4_+$, $t \geq 0$, $\bm{A} \in \mathbb{R}^{4 \times 4}$, $\bm{B} \in \mathbb{R}^4$, $u(t) \in \mathbb{R}$, $\bm{y}(t) \in \mathbb{\overline{R}}^l_+$, $\bm{C} \in \mathbb{R}^{l \times 4}$, and 
\begin{equation}\label{eq:AandB}
\bm{A} =
\begin{bmatrix}
-a_{41} & n_\text{q}a_{12} & n_\text{h}a_{13} & a_{14} \\
0 & \theta_1 & 0 & 0 \\
0 & 0 & \theta_2 & 0 \\
a_{41} & n_\text{q}a_{42} & n_\text{h}a_{43} & -a_{14} 
\end{bmatrix},
\quad
\bm{B} =
\begin{bmatrix}
-1 \\ 
0 \\ 
0 \\ 
0
\end{bmatrix},
\end{equation}
where $\theta_1 = -n_\text{q}(a_{12} + a_{42}) + a_{22}$ and $\theta_2 = -n_\text{h}(a_{13} + a_{43})+a_{33}$.
The dimension $l$ is the number of measured states. In Section \ref{sub:fullFeedback} we will cover the case of full-state feedback (i.e., $l = n = 4$), while Section \ref{sub:outputFeedback} will cover the case of output feedback (i.e., $l < n = 4$). The matrix $\bm{C}$ will take different forms in those two sections.

\subsection{Dynamics of Key Variables}
Throughout the paper, the variables depicting the dynamics of the climate and the natural resources are referred to as ``key variables''. The key variables we aim to regulate are: (1) the circularity of carbon dioxide and (2) the mean tropospheric temperature. As we will show, circularity quantifies two sustainability factors: first, the intensity of use of non-renewable resources, and second, the extent in which the system is close to the \emph{net-zero target} \cite{NZdefinition} requested by organizations such as the European Commission \cite{EU-netZero} and the United Nations \cite{UN-netZero}. In parallel, the mean tropospheric temperature is the main variable affecting climate change; the maximum temperature increase above the pre-industrial era was set to 1.5 $\tccentigrade$ in the Paris Agreement when it entered into effect on 4 November 2016 \cite{ParisAgreement}.

\subsubsection{Circularity of $\text{CO}_2$}
The dynamics of the \emph{instantaneous circularity}, namely, $\lambda(t)$, is given after the following definition.
\begin{definition}[\hspace{1sp}\cite{zocco2025cirl}]\label{def:unsMassFlow}
A mass or flow is \emph{finite-time sustainable} if either it exits a nonrenewable reservoir or it enters a landfill, an incinerator, or the natural environment as a pollutant.
\end{definition}

Thus, the \emph{instantaneous circularity} is given by \cite{zocco2025cirl}
\begin{equation}\label{eq:circularity}
\lambda(\mathcal{N}; t) = - \left(\Delta \mu_\text{fb} \overline{m}_{\textup{f}\textup{b}}(t) + \dot{\overline{m}}_{\textup{f}\textup{c}}(t)\right),
\end{equation} 
where $\mathcal{N}$ is a thermodynamical material network \cite{zocco2025cirl,zocco2023thermodynamical,zocco2025synchronized,zocco2025circular,zocco2022circularity}, $m_{\text{f}\text{b}}(t)$ is the \emph{net} finite-time sustainable mass transported in batches (e.g., solids transported on trucks), and $\dot{m}_{\text{f}\text{c}}(t)$ is the \emph{net} finite-time sustainable flow transported continuously (e.g., fluids transported through pipes). The quantity $\overline{m}_\text{fb}(t) = c_\text{fb} m_\text{fb}(t)$ is the \emph{weighted} $m_\text{fb}(t)$, where $c_\text{fb} \in (0, 1]$ is the criticality coefficient of the material indicating to what extent the material is critical (e.g., high value for a critical raw material \cite{CRMs-EU}), and $\mu_\text{fb} \geq 1$ takes into account how many batches of finite-time sustainable materials are disposed of while they are still \emph{functional} (i.e., disposing of faulty parts is better than disposing of working parts). Similarly, $\dot{\overline{m}}_\text{fc}(t) = c_\text{fc} \dot{m}_\text{fc}(t)$ is the \emph{weighted} $\dot{m}_\text{fc}(t)$, where $c_\text{fc} \in (0, 1]$ is the criticality coefficient of the finite-time sustainable material transported continuously. Finally, $\Delta > 0$ is a constant with the unit of a frequency introduced as a multiplying factor in order to convert the mass $\overline{m}_{\text{f}\text{b}}(t)$ into a flow, and hence, makes the sum with $\dot{\overline{m}}_{\text{f}\text{c}}(t)$ physically consistent. The choice of $\Delta$ is arbitrary, but its value must be kept the same for any calculation of $\lambda(\mathcal{N}; t)$ in order to make meaningful comparisons. Thus, $\lambda(\mathcal{N}; t) \in (-\infty, +\infty)$ and the circularity problem can be stated as
\begin{equation}
\mathcal{N}^* = \arg \max \,\, \lambda(\mathcal{N}; t). 
\end{equation}

For $\mathcal{N}_\text{c}$ shown in Fig. \ref{fig:compDiagraph}, note that $m_\text{fb}(t) = 0$. For $c_\text{fc}$ = 1, due to the high criticality of tropospheric $\text{CO}_2$, the instantaneous circularity (\ref{eq:circularity}) yields
\begin{equation}\label{eq:lambdaOfThisNet}
\lambda_\text{c}(\mathcal{N}_\text{c}; t) = -\dot{m}_\text{fc}(t) = - (\phi_1(t) + \phi_\text{nz}(t)),  
\end{equation}
where   
\begin{equation}
\phi_1(t) = a_{22} x_2(t) + a_{33} x_3(t),
\end{equation}
and
\begin{equation}
\begin{aligned}
\phi_\text{nz}(t) = n_\text{q} a_{12} x_2(t) + n_\text{q} a_{42} x_2(t) \\ + n_\text{h} a_{13} x_3(t) + n_\text{h} a_{43} x_3(t) - u(t).
\end{aligned}
\end{equation}
Note that, consistent with Definition \ref{def:unsMassFlow}, $\phi_1(t)$ takes into account the material flows extracted from non-renewable reserves (i.e., $c^6_{6,6}$), whereas $\phi_\text{nz}(t)$ is the difference between the mass flow rate of $\text{CO}_2$ emitted to the atmosphere (which is considered a pollutant because it is well beyond the natural atmospheric concentrations) and the mass flow rate of $\text{CO}_2$ removed by the carbon-capture system, i.e., $u(t)$. Specifically, $\phi_\text{nz}(t)$ corresponds to the net-zero target and it is met when
\begin{equation}\label{eq:phiNZ}
\phi_\text{nz}(t) = 0 \Rightarrow  x_1(t) + x_4(t) = \gamma,
\end{equation}
where $\gamma \in \overline{\mathbb{R}}_+$ is a constant. The constant $\gamma$ is the mass of $\text{CO}_2$ in the considered part of the troposphere, i.e., in the compartments $c^1_{1,1}$ and $c^4_{4,4}$ together, that, for a net-zero condition, have the total input flow equal to the total output flow (i.e., $u(t)$). Note that the exchange of carbon inside the troposphere, i.e., between $c^1_{1,1}$ and $c^4_{4,4}$, is not constrained in the net-zero target, and hence, in general $x_1(t)$ and $x_4(t)$ can vary over time even if condition (\ref{eq:phiNZ}) is reached.

\subsubsection{Mean Tropospheric Temperature}
The $\text{CO}_2$ concentration affects the average tropospheric temperature of the Earth. A simple model of the average tropospheric temperature consisting of a first-order nonlinear ordinary differential equation is given by \cite{elsherif2025climate}
\begin{equation}\label{eq:TemperatureDynamics}
\dot{T}(t) = \frac{1}{C} \left[ S(1 - \alpha) - \epsilon \sigma T(t)^4\right], 
\end{equation}
where $T(t)$ [K] is the mean tropospheric temperature, $C$ [J$\text{m}^{-2}\text{K}^{-1}$] is the effective heat capacity of the Earth, which quantifies the system sensitivity to energy variations, $S$ [J$\text{s}^{-1}\text{m}^{-2}$] is the average incoming solar energy per unit area, $\alpha \in [0, 1]$ models the fraction of solar energy reflected by the Earth (i.e., the albedo effect), $\sigma = 5.67 \times 10^{-8} \text{ J}\text{s}^{-1}\text{m}^{-2}\text{K}^{-4}$ is the Stefan-Boltzmann's constant, and $\epsilon \in [0, 1]$ is the emissivity of Earth, which quantifies how much energy is sent back to space and it is reduced by an increase of the concentration of $\text{CO}_2$ in the troposphere. Thus, the variation of the tropospheric temperature of Earth (\ref{eq:TemperatureDynamics}) is affected by two contrasting contributions: the amount of energy entering the troposphere and not reflected back (i.e., the term multiplied by $S$), and by the outgoing longwave radiation (i.e., $\epsilon \sigma T(t)^4$).

To take into account in (\ref{eq:TemperatureDynamics}) of the $\text{CO}_2$ concentration, which is regulated by our control scheme, the emissivity is modeled as  
\begin{equation}
\epsilon(t) = \epsilon_\text{n} - d_\text{c}(t),
\end{equation}
where
\begin{equation}
d_\text{c}(t) = \eta \frac{x_1(t)}{V_1},
\end{equation}
$\epsilon_\text{n}$ is the normal emissivity of a urban area, $d_\text{c}(t)$ is the disturbance to the emissivity induced by the carbon dioxide concentration in the target urban area $c^1_{1,1}$, $\eta = 2.36 \times 10^{-4} \text{ km}^3/\text{t}$ is a multiplying factor added to render $d_\text{c}(t)$ nondimensional while maintaining $\epsilon(t)$ within the usual range 0.6-1.0 \cite{NASAemissivity}, and $V_1$ is the volume of compartment $c^1_{1,1}$.

\subsection{Full-State Feedback}\label{sub:fullFeedback}
The control law leveraged in this paper for full-state feedback was proposed in \cite{nersesov2004optimal} and applied for optimal drug delivery. The analogy between drug delivery and carbon control for climate stability is that the states of the compartments must remain in the nonnegative orthant of the state space since they are masses, and hence, negative values are physically inconsistent. 

Given that the control input in our network system (Fig. \ref{fig:compDiagraph} and equations (\ref{eq:AxBu})-(\ref{eq:AandB})) is not constrained to be an injection for the regulated compartment (i.e., $c^1_{1,1}$), the controller in \cite[Section~6]{nersesov2004optimal} is not applicable. Thus, we implement the optimal non-zero set-point regulator covered in \cite[Section~4]{nersesov2004optimal} and we set $\bm{C}$ = $\bm{I}_n$, where $\bm{I}_n \in \mathbb{R}^{n \times n}$ denotes the $n \times n$ identity matrix. In this case, the control input takes the form
\begin{equation}\label{eq:controller}
u(t) = - \bm{K} (\bm{x}(t) -\bm{x}_\text{e}) + v_\text{e},
\end{equation}
where $\bm{K} \in \mathbb{R}^{1 \times n}$ is the control gain, $v_\text{e} \in \mathbb{R}$, $\bm{x}_\text{e} \in \overline{\mathbb{R}}^n_+$ is an equilibrium for (\ref{eq:AxBu}), and $\bm{x}_\text{e} = [x_{1\text{e}},x_{2\text{e}},x_{3\text{e}},x_{4\text{e}}]^\top$. Thus, the closed-loop system becomes
\begin{equation}\label{eq:closedLoopSys}
\dot{\bm{x}}(t) = \hat{\bm{A}}\left(\bm{x}(t) - x_\text{e}\right), \,\, \bm{x}(0) = \bm{x}_0, \,\, t \geq 0,
\end{equation} 
where $\hat{\bm{A}} = \bm{A} - \bm{BK}$.
If
\begin{itemize}
    \item{the pair ($\bm{A},\bm{B}$) is stabilizable,}
    \item{[set-point condition] there exists $v_\text{e} \in \mathbb{R}$ such that $0 = \bm{A}\bm{x}_\text{e} + \bm{B} v_\text{e}$, where $\bm{x}_\text{e} \in \overline{\mathbb{R}}^n_+$ is the desired set-point,}
    \item{the pair ($\bm{A}, \bm{R}_1$) is observable, with $\bm{R}_1 \in \mathbb{R}^{n \times n}$ and $\bm{R}_1 \geq 0$,}
\end{itemize}

then 
\begin{itemize}
\item{there exists $\bm{P} > 0$ such that 
\begin{equation}\label{eq:ARE}
    0 = \bm{A}^\top \bm{P} + \bm{PA} + \bm{R}_1 - \bm{PSP},
    \end{equation}
    where $\bm{S} = r_2^{-1}\bm{B}\bm{B}^\top$ and $r_2 > 0$,}
\item{and controller $\bm{K} \in \mathbb{R}^{1 \times n}$ (a) guarantees that the closed-loop system (\ref{eq:closedLoopSys}) is asymptotically stable and (b) minimizes the performance functional
\begin{equation}\label{eq:performanceFunctional}
\begin{aligned}[b]
J(\bm{x}_0, u) & = \frac{1}{2}\int_0^\infty \left[(\bm{x}(t) - \bm{x}_\text{e})^\top \bm{R}_1 (\bm{x}(t) - \bm{x}_\text{e}) \right. \\ & \left. + r_2 (u(t) - v_\text{e})^2 \right] \text{d}t.
\end{aligned}
\end{equation}
}
\end{itemize}
Furthermore, $\bm{x}(t) \in \mathcal{D}_\text{A} \subseteq \overline{\mathbb{R}}^{n}_+$ for $t \geq 0$ and every $\bm{x}_0 \in \mathcal{D}_\text{A}$, where $\mathcal{D}_\text{A} \subseteq \overline{\mathbb{R}}^{n}_+$ is a subset of the domain of attraction of (\ref{eq:closedLoopSys}) given in \cite[Theorem 4.1]{nersesov2004optimal}. If $\bm{x}_0 \in \mathcal{D}_{\text{A}}$, then $\bm{K} \in \mathbb{R}^{1 \times n}$ guarantees that $\bm{x}(t) \geq \geq 0$, that is, the state trajectories remain in the non-negative orthant of $\mathbb{R}^n$ at all time, which preserves the fundamental property of non-negative states for the closed-loop system (\ref{eq:closedLoopSys}). 

Hence, the procedure to design the controller is as follows.
\begin{enumerate}
\item{Verify that the pair ($\bm{A}$, $\bm{B}$) is stabilizable.}
\item{[Set-point condition] Find $v_\text{e} \in \mathbb{R}$ such that $0 = \bm{A}\bm{x}_\text{e} + \bm{B} v_\text{e}$, where $\bm{x}_\text{e}$ is the desired set-point.}
\item{Choose $\bm{R}_1 \geq 0$ such that the pair $(\bm{A}, \bm{R}_1)$ is observable.}
\item{Choose $r_2 > 0$ and solve the algebraic Riccati equation (\ref{eq:ARE}).}
\item{Compute $\bm{K} = (1/r_2) \bm{B}^\top \bm{P}$.}
\item{Compute the translated initial conditions as $\tilde{\bm{x}}_0 = \bm{x}_0 - \bm{x}_\text{e}$ and choose the initial conditions $\bm{x}_0$ such that either $\bm{x}_0 \in \mathcal{D}_A$, with $\mathcal{D}_A$ calculated as in \cite[Theorem 4.1]{nersesov2004optimal}, or the closed-loop state $\bm{x}(t) \geq \geq 0$, $t \geq 0$, which can be verified a posteriori.}
\item{Implement the closed-loop system as 
\begin{equation}\label{eq:transClosedLoop}
\dot{\tilde{\bm{x}}}(t) = \hat{\bm{A}}\tilde{\bm{x}}(t), \,\, \tilde{\bm{x}}(0) = \tilde{\bm{x}}_0, \,\, t \geq 0,
\end{equation}
where $\tilde{\bm{x}}(t) = \bm{x}(t) - \bm{x}_\text{e}$ is the translated state so that the origin $\bm{0}$ corresponds to the desired set-point $\bm{x}_\text{e}$.} 
\item{Recover the original state as $\bm{x}(t) = \tilde{\bm{x}}(t) + \bm{x}_\text{e}$ and the original initial conditions as $\bm{x}_0 = \tilde{\bm{x}}_0 + \bm{x}_\text{e}$.}
\item{Verify that $\lim_{t \to \infty} u(t) = v_e$ by computing $u(t)$ given by (\ref{eq:controller}); $u(t) > 0, t \geq 0$, indicates that the direction of $u(t)$ is the one considered in (\ref{eq:AandB}).}
\item{Verify that the closed-loop state $\bm{x}(t) \geq \geq 0$; if not, $\bm{x}_\text{e}$ and/or $\bm{x}_0$ must be changed.}
\end{enumerate}

The set-point condition 2) fixes at which state $\bm{x}_\text{e}$ the closed-loop system is desired to converge. After a transitory phase, the closed-loop system states converge to the equilibrium $\bm{x}_\text{e}$ while the control input $u(t)$ converges to the constant value $v_e$. Note that the origin of the dynamical system (\ref{eq:dynamicEq1})-(\ref{eq:dynamicEq4}), i.e., $\bm{x}(t) \equiv \bm{0}$, is not a desired set-point for our carbon control problem. Specifically, we need to choose $x_{i\text{e}}|_{i \in \{1,4\}} \neq 0$ to conform with real-world conditions, i.e., the masses of $\text{CO}_2$ in the troposphere compartments ($c^1_{1,1}$ and $c^4_{4,4}$) are nonzero. Thus, imposing the set-point condition to the dynamical system (\ref{eq:dynamicEq1})-(\ref{eq:dynamicEq4}) with $x_{i\text{e}}|_{i \in \{1,4\}} \neq 0$ yields
\begin{gather}\label{eq:equilibriumCond1}
v_e = -a_{41} x_{1\text{e}} + a_{14} x_{4\text{e}} + n_{\text{h}} a_{13} x_{3\text{e}} + n_\text{q} a_{12} x_{2\text{e}}, 
\\ \label{eq:equilibriumCond2}
0 = \left(a_{22} - n_\text{q} a_{12} - n_\text{q} a_{42}\right) x_{2\text{e}}, 
\\ \label{eq:equilibriumCond3}
0 = \left(a_{33} - n_\text{h} a_{13} - n_\text{h} a_{43}\right) x_{3\text{e}},
\\ \label{eq:equilibriumCond4}
0 = a_{41} x_{1\text{e}} - a_{14} x_{4\text{e}} + n_\text{q} a_{42} x_{2\text{e}} + n_\text{h} a_{43} x_{3\text{e}}.
\end{gather}
Conditions (\ref{eq:equilibriumCond2}) and (\ref{eq:equilibriumCond3}) require that $x_{i\text{e}}|_{i \in \{2,3\}} = 0$ (i.e., no mass of fuel and gas in vehicles and houses at equilibrium); if $x_{i\text{e}}|_{i \in \{2,3\}} \neq 0$, conditions (\ref{eq:equilibriumCond2}) and (\ref{eq:equilibriumCond3}) yield $\theta_1 = \theta_2 = 0$, which renders the pair $(\bm{A},\bm{B})$ unstabilizable. With $x_{i\text{e}}|_{i \in \{2,3\}} = 0$, equations (\ref{eq:equilibriumCond1})-(\ref{eq:equilibriumCond4}) become
\begin{gather}\label{eq:setPoint1}
v_\text{e} = 0,  
\\ \label{eq:setPoint2}
a_{22} \neq n_\text{q}(a_{12}+a_{42}), \quad \text{(i.e., } \theta_1 \neq 0) 
\\ \label{eq:setPoint3}
a_{33} \neq n_\text{h}(a_{13}+a_{43}), \quad \text{(i.e., } \theta_2 \neq 0)
\\ \label{eq:setPoint4} 
x_{4\text{e}} = \frac{a_{41}}{a_{14}} x_{1\text{e}},
\end{gather}
with $\bm{x}_\text{e} = [x_{1\text{e}}, 0, 0, x_{4\text{e}}]^\top$.

\subsection{Output Feedback}\label{sub:outputFeedback}
In practice, it is often unfeasible to measure the full state of a system as required for a full-state feedback controller. When the full state is not accessible for feedback, the LQR design problem is referred to as ``static output feedback''. In this section, we cover the case in which only the mass of $\text{CO}_2$ in the target area is known, i.e., in the compartment $c^1_{1,1}$. Thus,
\begin{equation}
\bm{y}(t) = \bm{C} \bm{x}(t), \quad \bm{C} = [1, 0, 0, 0],
\end{equation}
while the control input is given by
\begin{equation}\label{eq:controller-OutputFdbk}
u(t) = -\bm{K}\left(\bm{y}(t) - \bm{C}\bm{x}_\text{e} \right) + v_\text{e},  
\end{equation}
where $\bm{K} \in \mathbb{R}^{1 \times l}$ is the controller gain. When $\bm{C} \neq \bm{I}_n$, finding an optimal stabilizing feedback gain $\bm{K}$ is more complex and requires the solution of a set of two coupled equations; one modified Riccati equation and one modified Lyapunov equation \cite{nersesov2004optimal}. Here, we use the algorithms developed by Ilka and Murgovski \cite{ilka2022novel} to solve the static output feedback control problem. The procedure is similar to the LQR design and can be summarized as follows with the conditions for stabilizability detailed in \cite[Theorem 1]{ilka2022novel}.
\begin{enumerate}
\item{Verify that the pair ($\bm{A}, \bm{B}$) is stabilizable (e.g., via the PBH test).}
\item{Verify that the pair ($\bm{A}, \bm{C}$) is detectable (e.g., via the PBH test).}
\item{Verify that the weighting matrices in the performance functional (\ref{eq:performanceFunctional}) satisfy
\begin{equation}\label{eq:conditionOfWeightMatrices}
\begin{bmatrix}
\bm{R}_1  &  \bm{0}\\
\bm{0}^\top  & r_2\\
\end{bmatrix} \geq 0, \quad r_2 > 0,
\end{equation}
with $\bm{\bm{R}}_1 \in \mathbb{R}^{n \times n}$ and $r_2 \in \mathbb{R}$ such that $\bm{\bm{R}}_1 \geq 0$ and $r_2 >0$, and where $\bm{0} \in \mathbb{R}^{n \times 1}$ is a vector of zeros.
}
\item{Run the algorithm given in \cite{ilka2022novel} to find a feedback gain $\bm{K} \in \mathbb{R}^{1 \times l}$ and a solution $\bm{P} \in \mathbb{R}^{n \times n}$ to the algebraic Riccati equation
\begin{equation}\label{eq:AREoutputFdbk}
\begin{aligned}
\bm{A}^\top\bm{P} + \bm{P}\bm{A} + \bm{R}_1 + r_2\bm{G}^\top \bm{G} \\ - r^{-1}_2 \bm{PB} \bm{B}^\top \bm{P} = \bm{0},
\end{aligned}
\end{equation}
where $\bm{G} \in \mathbb{R}^{1 \times n}$ is computed as  
\begin{equation}\label{eq:matrixG}
\bm{G} = \bm{KC} - r^{-1}_2 \bm{B}^\top \bm{P}.
\end{equation}
}
\item{Compute $\bm{G}\in \mathbb{R}^{1 \times n}$ as in (\ref{eq:matrixG}).}
\item{Check that (\ref{eq:AREoutputFdbk}) holds.}
\item{Implement the closed-loop system as in (\ref{eq:transClosedLoop}) with $\hat{\bm{A}} = \bm{A} - \bm{BKC}$.} 
\item{Recover the original state as $\bm{x}(t) = \tilde{\bm{x}}(t) + \bm{x}_\text{e}$ and the original initial conditions as $\bm{x}_0 = \tilde{\bm{x}}_0 + \bm{x}_\text{e}$.}
\item{Verify that $\lim_{t \to \infty} u(t) = v_e$ by computing $u(t), t \geq 0$, as in (\ref{eq:controller-OutputFdbk}) with $\bm{y}(t) = \bm{C}(\tilde{\bm{x}}(t) + \bm{x}_\text{e})$.}
\item{Verify that the closed-loop state $\bm{x}(t) \geq \geq 0, t \geq 0$; if not, $\bm{x}_\text{e}$ and/or $\bm{x}_0$ must be changed.}
\end{enumerate}

\section{Numerical Studies: Full-State Feedback}\label{sec:NumStudies-fullStateFdbk}
This section is divided into two parts. In the first part, we consider arbitrary rate constants in the linear system (\ref{eq:AandB}) for testing the control scheme. In the second part, we move to consider realistic parameters and address the control of tropospheric carbon (\ref{eq:TemperatureDynamics}) and the system circularity (\ref{eq:lambdaOfThisNet}). For the full-state feedback, we leveraged the Python library of Fuller \emph{at al.} \cite{fuller2021python}.

\subsection{Arbitrary Rate Constants}\label{sub:FullStateFdbk-RandomRates}
\begin{figure}
\centering
\subfloat{
  \includegraphics[width=0.23\textwidth]{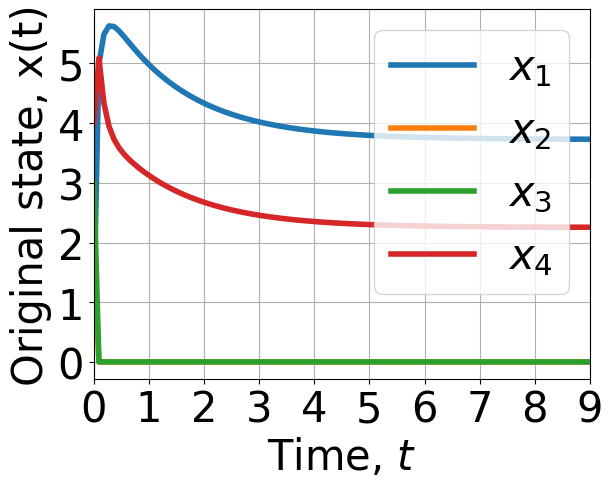}
  }
\subfloat{
  \includegraphics[width=0.23\textwidth]{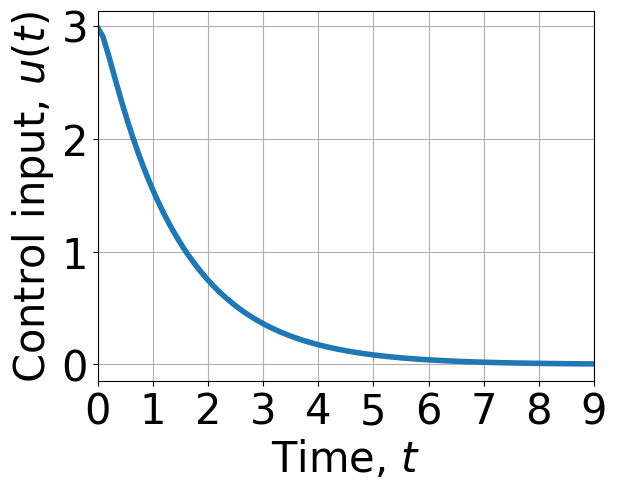}
  }
\caption{Closed-loop state and control input with \emph{full-state feedback} and arbitrary rate constants $a_{k,s}$.}
\label{fig:fullState-Rdm}
\end{figure}

We begin the numerical studies by considering the following state matrix with arbitrary rate constants $a_{k,s}$, with $n_\text{q}$ = 15, $n_\text{h}$ = 7, and 
\begin{equation}\label{eq:AwithRandomRates}
\bm{A} =
\begin{bmatrix}
-3.41 & 95.82 & 39.62 & 5.63 \\
0 & -128.85 & 0 & 0 \\
0 & 0 & -95.20 & 0 \\
3.41 & 38.57 & 63.91 & -5.63 
\end{bmatrix}.
\end{equation}
Thus, $\theta_1$ = -128.85 and $\theta_2$ = -95.20. The matrix $\bm{A}$ has eigenvalues 0, -9.04, -128.85, and -95.20, and using the Popov–Belevitch–Hautus (PBH) test \cite{williams2007linear}, the pair ($\bm{A}, \bm{B}$) is stabilizable. The set-point is also set by arbitrarily choosing $\bm{x}_\text{1e}$ = 3.73 and imposing (\ref{eq:setPoint1})-(\ref{eq:setPoint4}); thus, $\bm{x}_\text{e} = [3.73, 0, 0, 2.25]^\top$. We set $\bm{R}_1$ = $\bm{I}_4$, and hence, the pair ($\bm{A}, \bm{R}_1$) is observable. The solution to the algebraic Riccati equation (\ref{eq:ARE}) with $r_2$ = 1 yields a positive-definite $\bm{P}$, which is then used to compute the feedback gain $\bm{K}$ = [-0.75, -0.76, -0.77, -0.69]. With testing intial conditions $\bm{x}_0 = [1,2,3,4]^\top$, the closed-loop original states and control input are shown in Fig. \ref{fig:fullState-Rdm}. The states converge to the desired set-point $\bm{x}_\text{e}$ with the equilibrium input $v_\text{e} = 0$. Moreover, $\bm{x}(t) \geq \geq 0, t \geq 0$, which is physically consistent since the states represent masses.

\subsection{Case of Tropospheric $\text{CO}_2$}\label{sub:FullState-CO2case}
We now consider realistic system initial conditions and set-point to address the control of the concentration of tropospheric carbon. Specifically, we consider the current state of the atmosphere as our initial condition and the pre-industrial condition as the set-point state. In other words, the controller is expected to bring the concentration of $\text{CO}_2$ to the pre-industrial era. 

Let $x_{i0}$ be the initial condition of the $i$-th state. We consider a target urban area with volume $V_1 = 1.18 \text{ km}^3$, which results from the product between the surface of the Italian city of Siena (approximately 118 $\text{km}^2$) and a hight of 0.01 km. Given that the current concentration of carbon dioxide in the atmosphere is 431 ppm \cite{NASAcarbon}, the current mass of carbon inside the volume $V_1$ is 915.4 tonnes; thus, $x_{10}$ = 915.4 t. The surrounding area $c^4_{4,4}$ has a volume twice the one of $c^1_{1,1}$, and hence, $x_{40}$ =  1,830.8 t. To set the initial mass of fuel (diesel) in $c^2_{2,2}$, i.e., $x_{20}$, we consider $n_\text{q} = 5,000$ vehicles with a gas tank of 50 liters each. This yields 42 kg of diesel per vehicle, and hence, $x_{20}$ = 210 t. The initial mass of gas inside pipelines and in one house is set as 50 kg, and hence, with the number of houses set as $n_\text{h} = 10,000$, we have that $x_{30}$ = 500 t. 

\begin{figure}
\centering
\subfloat{
  \includegraphics[width=0.23\textwidth]{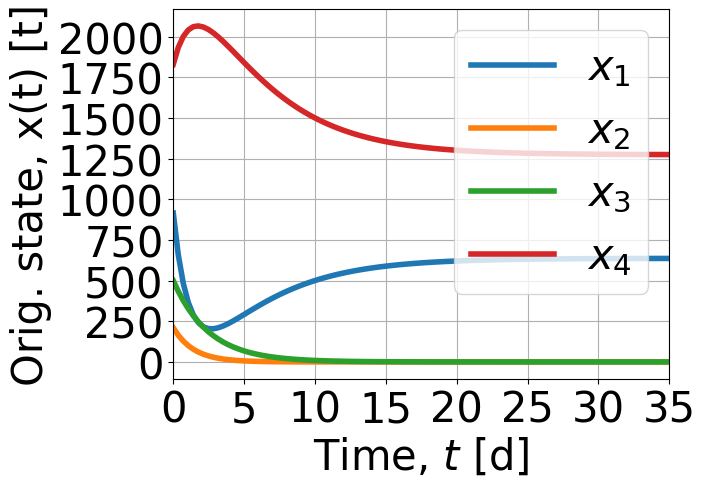}
  }
\subfloat{
  \includegraphics[width=0.23\textwidth]{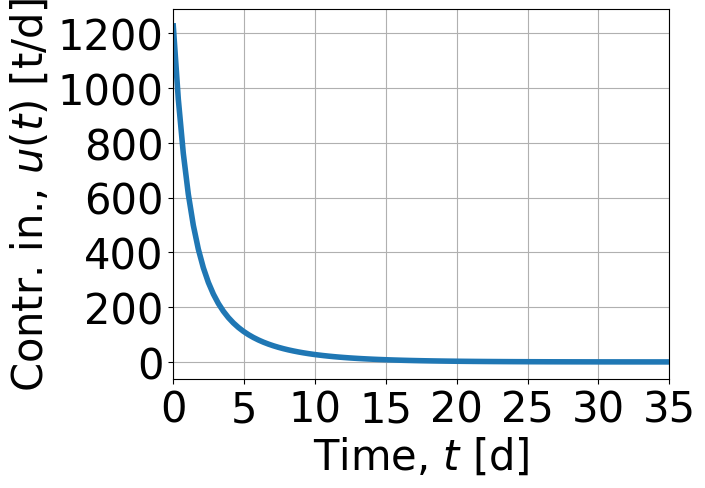}
  }
\caption{Closed-loop state and control input with \emph{full-state feedback} for tropospheric $\text{CO}_2$ control.}
\label{fig:FullState_CO2}
\end{figure}

Since the set-point condition corresponds to the pre-industrial era in which the concentration of $\text{CO}_2$ was approximately 300 ppm, we set $x_{1\text{e}}$ = 637.2 t. Assuming the rate coefficients in Table \ref{tab:rateConstantsCO2}, we have that $x_{4\text{e}} = (a_{41}/a_{14})x_{1\text{e}}$ = 1,274.4 t (imposition of (\ref{eq:setPoint4})). 
\begin{table}
\centering
\caption{Rate constants used in the case of tropospheric $\text{CO}_2$.}
\label{tab:rateConstantsCO2}
\begin{tabular}{cccccccc} 
\hline
$a_{41}$ & $a_{12}$ & $a_{13}$ & $a_{14}$ & $a_{42}$ & $a_{22}$ & $a_{43}$ & $a_{33}$\\ [1.5ex] 
 
 0.2 & 0.5/$n_\text{q}$ & 0.5/$n_\text{h}$ & 0.1 & 0.5/$n_\text{q}$ & 0.3 & 0.5/$n_\text{h}$ & 0.6 \\ 
\hline
\end{tabular}
\end{table}
Thus, the system matrix (\ref{eq:AandB}) becomes
\begin{equation}\label{eq:AinCO2case}
\bm{A} =
\begin{bmatrix}
-0.2 & 0.5 & 0.5 & 0.1 \\
0.0 & -0.7 & 0.0 & 0.0 \\
0.0 & 0.0 & -0.4 & 0.0 \\
0.2 & 0.5 & 0.5 & -0.1 \\
\end{bmatrix},
\end{equation}
with $\theta_1$ = -0.7 and $\theta_2$ = -0.4. The eigenvalues of $\bm{A}$ are -0.3,  0,  -0.7, and -0.4, and the pair ($\bm{A}, \bm{B}$) is stabilizable. By choosing $\bm{R}_1$ = $\bm{I}_4$, the pair ($\bm{A}, \bm{R}_1$) is observable. With $r_2 = 1$, the optimal gain is $\bm{K}$ = [-0.94, -0.67, -0.93, -0.65]. The response of the closed-loop system to the initial conditions specified previously and the control input are shown in Fig. \ref{fig:FullState_CO2}. As can be seen, \emph{the concentration in the pre-industrial era $x_{1\text{e}}$ is reached in approximately 25 days}. Note that, since $x_{2\text{e}}$ = $x_{3\text{e}}$ = 0, all the mass flow rates converge to zero with the exception of the flows between $c^1_{1,1}$ and $c^4_{4,4}$, i.e., the tropospheric regions are a closed system at equilibrium. This result explains why the paper title specifies that the network is ``finite-time combustion-tolerant'': to reach the desired pre-industrial era condition inside the target area, i.e., to reach $x_{1\text{e}}$, this control scheme brings to zero the $\text{CO}_2$ emissions of the $n_\text{q}$ vehicles and the $n_\text{h}$ houses, i.e, the mass flow rates from $c^2_{2,2}$ and $c^3_{3,3}$ to $c^1_{1,1}$ and $c^4_{4,4}$. Specifically, \emph{this full-state control scheme tolerates the combustion taking place in 5,000 vehicles and in 10,000 house heating systems for approximately 6 days}, which is the time needed to bring $x_2(t)$ and $x_3(t)$ to zero.       

\begin{figure*}
\centering
\subfloat{
  \includegraphics[width=0.32\textwidth]{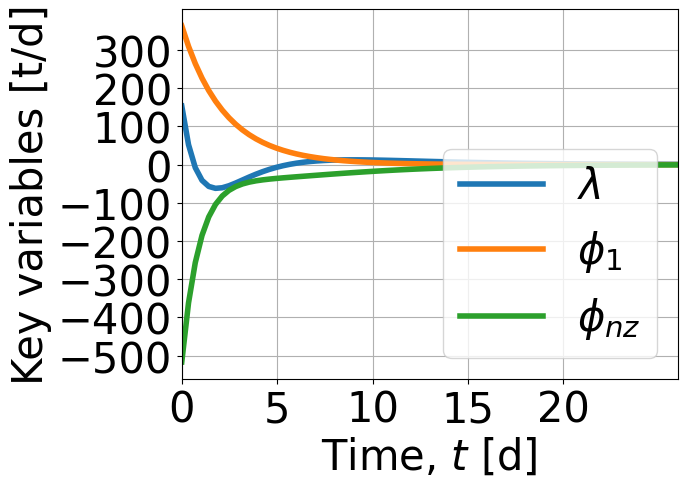}
  }
\subfloat{
  \includegraphics[width=0.27\textwidth]{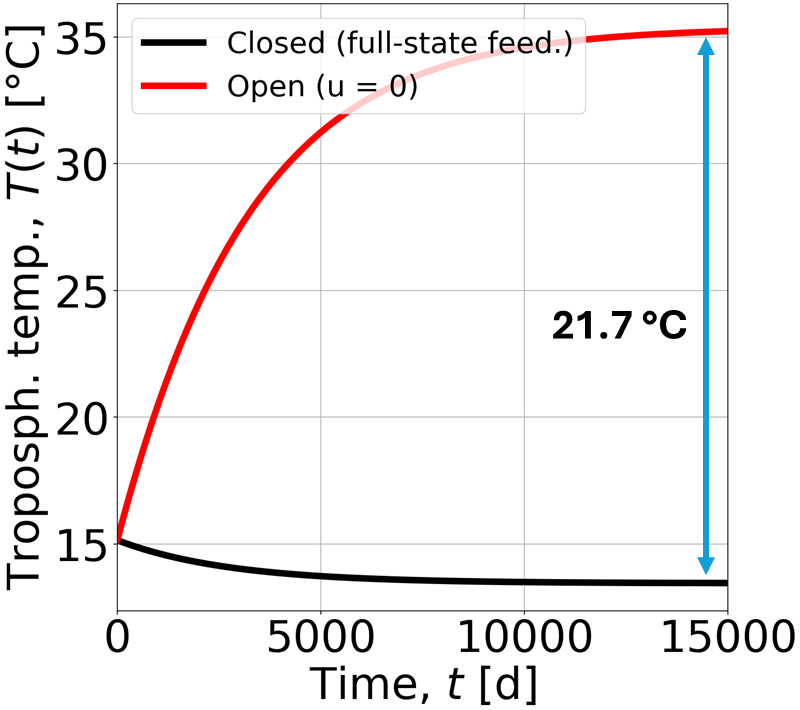}
  }
  \subfloat{
  \includegraphics[width=0.27\textwidth]{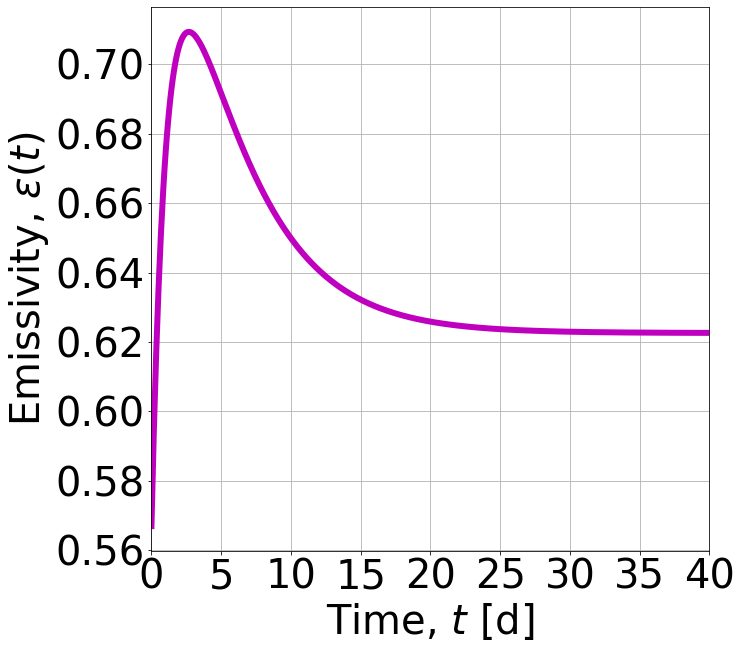}
  }
\caption{Key variables with \emph{full-state feedback}: (left) circularity, $\phi_1$, and $\phi_\text{nz}$; (middle) tropospheric temperature with open and closed loops; (right) emissivity.}
\label{fig:FullState_CO2_keyVariables}
\end{figure*}

The dynamics of the key variables are shown in Fig. \ref{fig:FullState_CO2_keyVariables}. Initially, circularity $\lambda(t)$ is positive because the extraction of carbon, i.e., $u(t)$, dominates over the other flows; at equilibrium, the extraction of finite natural resources $\phi_1(t)$ reaches zero, and so does the net-zero variable $\phi_{\text{nz}}(t)$. In particular, this controller meets the net-zero target (\ref{eq:phiNZ}) by nullifying  both the input flows and the output flows of the tropospheric areas $c^1_{1,1}$ and $c^4_{4,4}$. Fig. \ref{fig:FullState_CO2_keyVariables} also shows the dynamics of the tropospheric temperature, for which we set the current atmospheric temperature reported by NASA \cite{NASAtemperature} as the initial condition, i.e., $T(0)$ =  15.2 °C. The simulation is obtained for $C$ = 8 $\times$ 10$^8$ J K$^{-1}$ m$^{-2}$ \cite{elsherif2025climate}, $\alpha$ = 0.3 \cite{elsherif2025climate}, and $S$ = 1361/4 $\times$ 86400 J d$^{-1}$ m$^{-2}$ \cite{NASAirradiance}. The temperature has dynamics much slower than that of the states since its first derivative starts to reduce after approximately 5000 days, that is, 13.7 years; in the uncontrolled case (i.e., no carbon capture $u$ = 0), the temperature reaches the steady-state value of 35.2 °C in approximately 15000 days (41 years). In contrast, the steady-state value is 13.5 °C with full-state feedback control, which is 21.7 °C lower than the case of no carbon capture. Note that the steady-state temperature with closed loop (13.5 °C) is very close to the tropospheric temperature in the pre-industrial era (13.9 °C), which shows that the controller brings the temperature back to the pre-industrial era condition, as desired. Finally, Fig. \ref{fig:FullState_CO2_keyVariables} shows the emissivity $\epsilon(t)$, which is inversely proportional to $\bm{x}_{1}(t)$. 

\textbf{Remark:} In terms of the implementation of (\ref{eq:TemperatureDynamics}) in combination with the closed-loop system, it is worth pointing out that the former is nonlinear, whereas the latter is linear. Thus, we achieved the dynamics of $T(t)$ by numerically integrating the resulting system of five ordinary differential equations with the Python function \emph{scipy.integrate.odeint()}. Coherently with the states $x_i(t)|_{i = 1,2,3,4}$, we translated the temperature to the condition in the pre-industrial era, and hence, $\tilde{T}(t)$ = $T(t)$ - $T_\text{e}$, with $T_\text{e}$ = 13.9 °C. However, $T_\text{e}$ may not be an equilibrium for (\ref{eq:TemperatureDynamics}), which is not required by the control law because the regulator acts on the states $x_i(t)|_{i = 1,2,3,4}$, and not directly on $T(t)$. The controller affects $T(t)$ indirectly through $\epsilon(t)$, which is a function of $x_1(t)$.

\section{Numerical Studies: Output Feedback}\label{sec:NumStudies-outputFdbk}
This section presents numerical studies for the output feedback case with the controller design detailed in Section \ref{sub:outputFeedback}. As in Section \ref{sec:NumStudies-fullStateFdbk}, we first test the control with arbitrarily chosen rate constants $a_{k,s}$, and then, consider and discuss the case of rate constants set to tackle the control of the tropospheric carbon and the carbon circularity in a target area $c^1_{1,1}$. For the output feedback problem, we leveraged the MATLAB algorithm implementation made available by Ilka and Murgovski \cite{ilka2022novel}.

\subsection{Arbitrary Rate Constants}\label{sub:OutputFdbk-RandomRates}
\begin{figure*}
\centering
\subfloat{
  \includegraphics[width=0.3\textwidth]{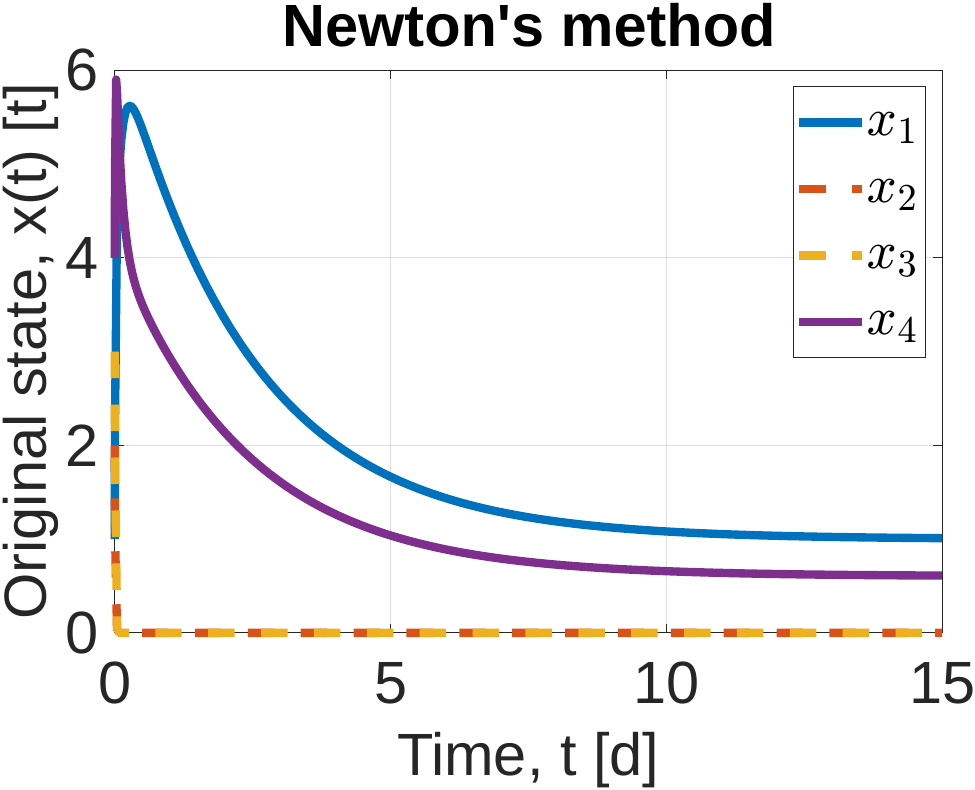}
  }
\subfloat{
  \includegraphics[width=0.3\textwidth]{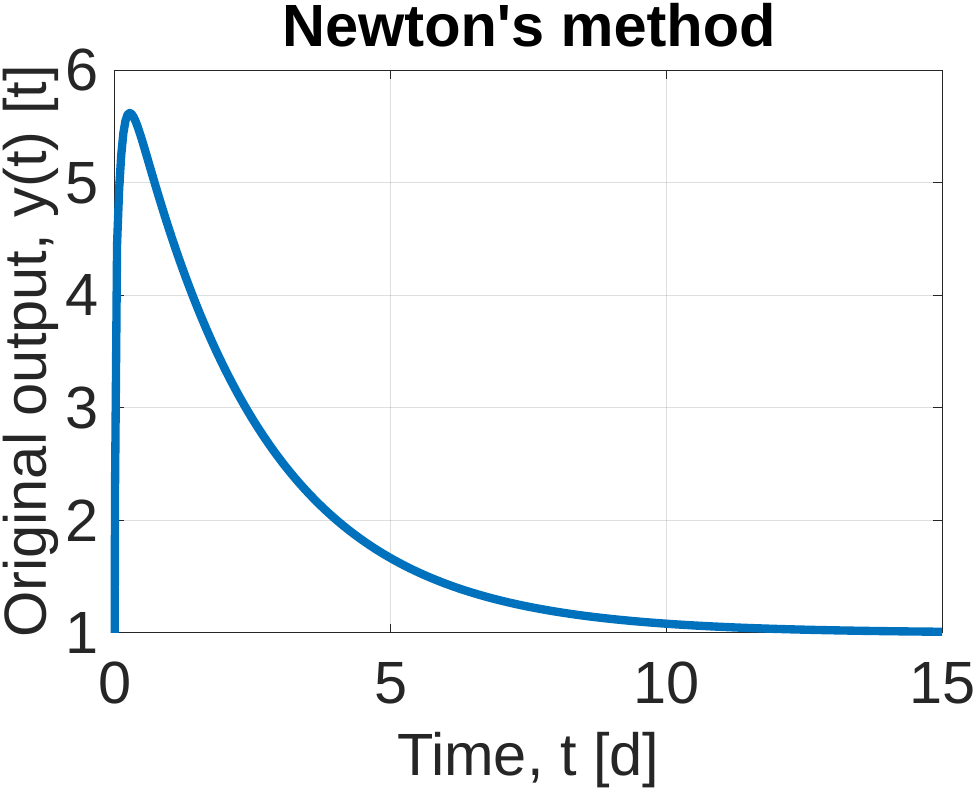}
  }
  \subfloat{
  \includegraphics[width=0.3\textwidth]{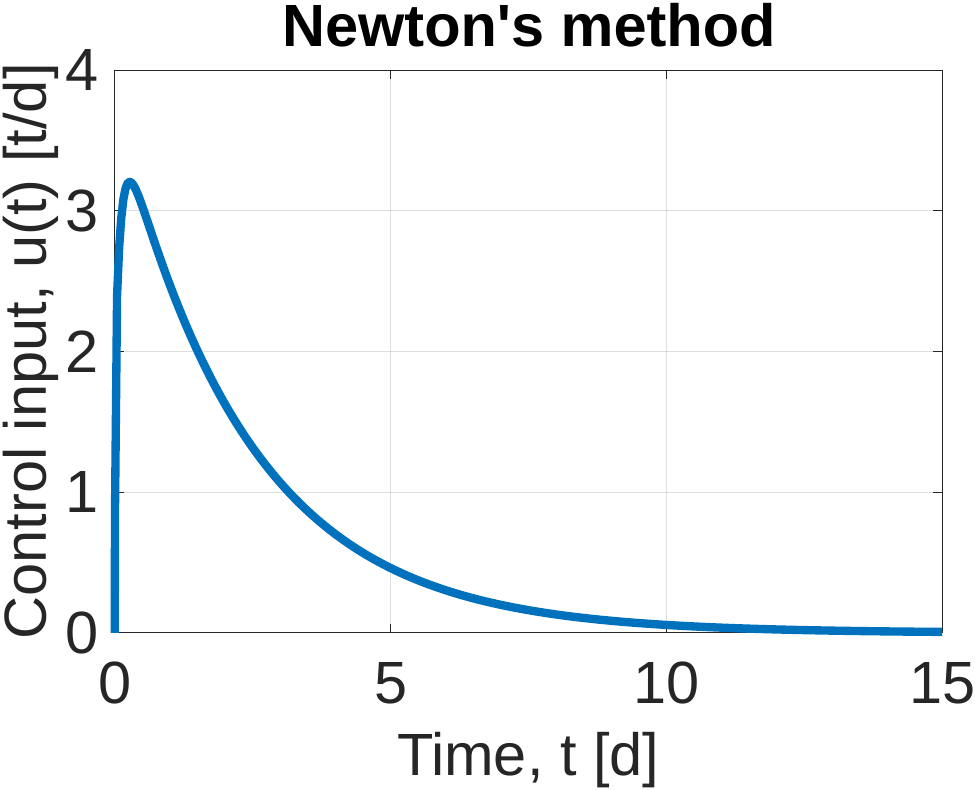}
  }
\caption{Closed-loop state, output, and control input with \emph{output-feedback} gain computed via the Newton's method \cite{ilka2022novel} in the case of arbitrary rate constants $a_{k,s}$.}
\label{fig:outputFdbk_rndRates}
\end{figure*}

We consider the same rate constants given in Section \ref{sub:FullStateFdbk-RandomRates}, and hence, the system matrix $\bm{A}$ is given by (\ref{eq:AwithRandomRates}). Thus, the pair ($\bm{A}, \bm{B}$) is stabilizable. The PBH test shows also that the pair ($\bm{A}, \bm{C}$) is detectable. As in \cite{ilka2022novel}, the weighting matrices have been chosen as $\bm{R}_1 = \bm{C}^\top \bm{C}$ and $r_2 = 1$. Thus, condition (\ref{eq:conditionOfWeightMatrices}) is satisfied. By running Newton's method and the modified Newton's method of \cite{ilka2022novel}, we obtained $\bm{K} = -0.694$ with both algorithms. The two methods returned slightly different $\bm{P}$ matrices, with the difference between corresponding entries of the two returned $\bm{P}$ matrices being smaller than 10$^{-7}$, and hence, both methods gave approximately the same $\bm{G}$ = [0, 0.696, 0.698, 0.618]. Furthermore, condition (\ref{eq:AREoutputFdbk}) was satisfied by both methods since the left-hand side of the equality contained entries smaller than 10$^{-7}$. The implementation of the closed-loop system (\ref{eq:closedLoopSys}) with $\hat{\bm{A}} = \bm{A} - \bm{BKC}$, testing initial conditions $\bm{x}_0$ = [1, 2, 3, 4]$^\top$, and set-point state $\bm{x}_\text{e}$ = [1, 0, 0, 0.605]$^\top$ gave the results shown in Fig. \ref{fig:outputFdbk_rndRates}. The controller stabilizes the translated state to the origin because its goal is to minimize the performance functional (\ref{eq:performanceFunctional}), while for the original state and control input we have that $\lim_{t \to \infty} \bm{x}(t) = \bm{x}_e$ and $\lim_{t \to \infty} u(t) = v_\text{e}$ (recall that $v_\text{e} = 0$ from equation (\ref{eq:setPoint1})). In addition, the closed-loop system state $\bm{x}(t) \ge \geq 0, t \geq 0$, which is physically consistent.       

\subsection{Case of Tropospheric $\text{CO}_2$}\label{sub:OutputFdbk-CO2case}
\begin{figure*}
\centering
\subfloat{
  \includegraphics[width=0.3\textwidth]{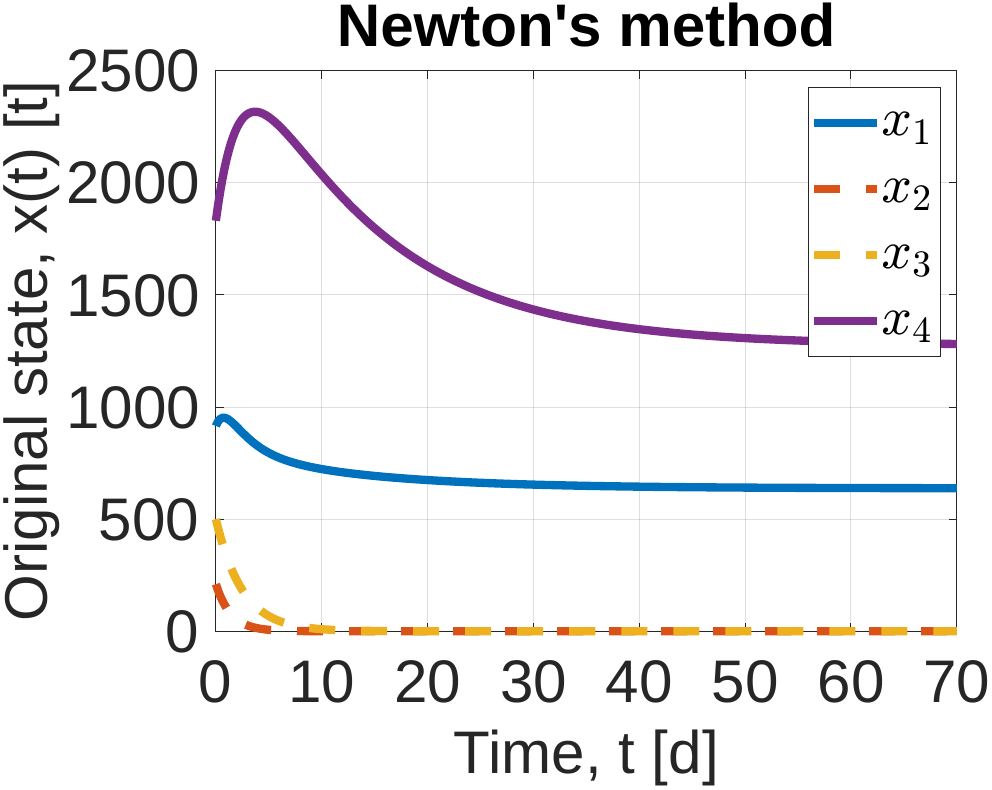}
  }
\subfloat{
  \includegraphics[width=0.3\textwidth]{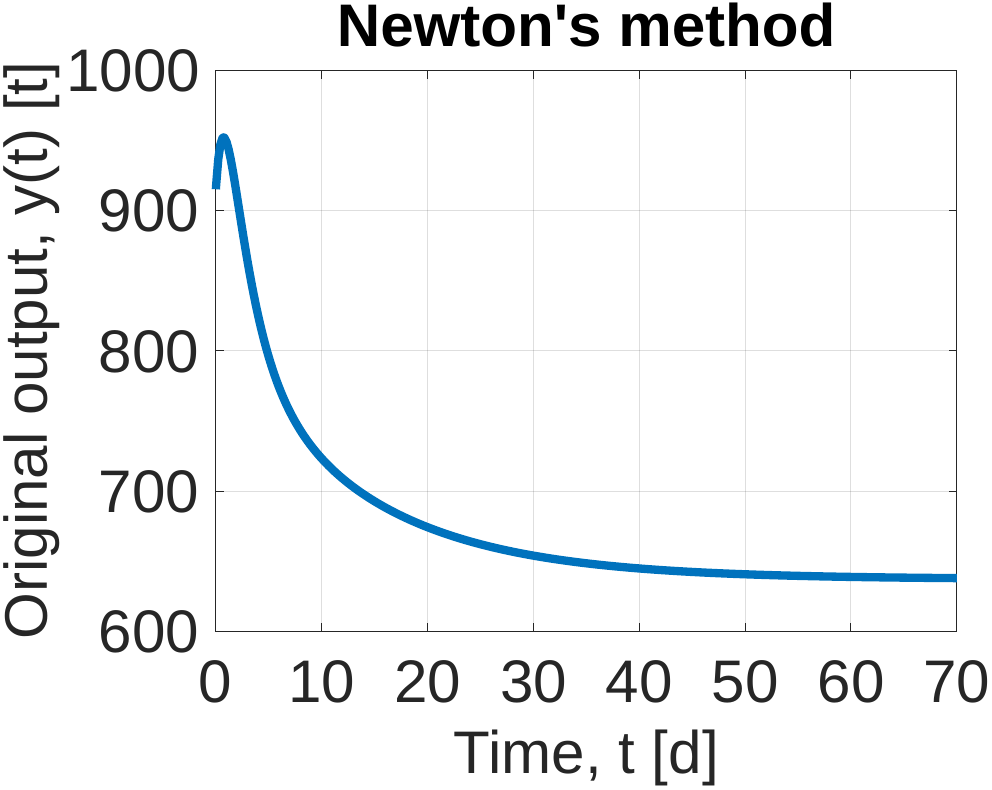}
  }
  \subfloat{
  \includegraphics[width=0.3\textwidth]{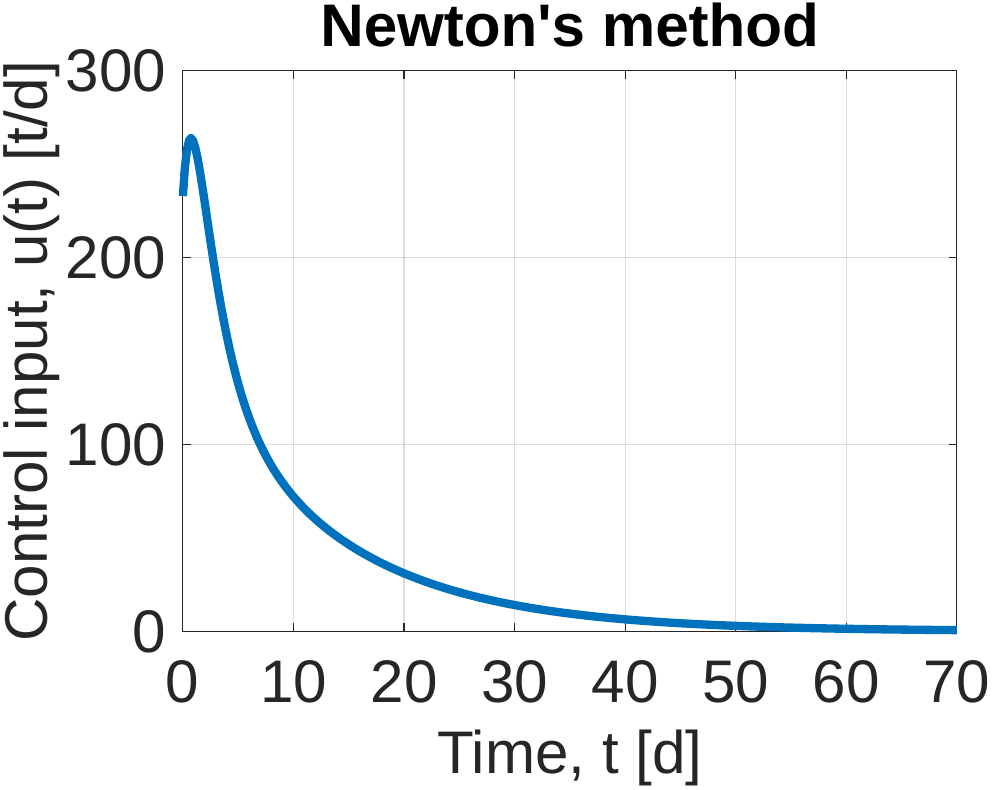}
  }
\caption{Closed-loop state, output, and control input with \emph{output-feedback} gain computed via the Newton's method \cite{ilka2022novel} for tropospheric $\text{CO}_2$ control.}
\label{fig:outputFdbk_CO2}
\end{figure*}

This section addresses the control of tropospheric carbon as in Section \ref{sub:FullState-CO2case} with output-feedback control. Thus, the system matrix is given by (\ref{eq:AinCO2case}) and the control design procedure follows that of Section \ref{sub:OutputFdbk-RandomRates} with $\bm{R}_1 = \bm{C}^\top \bm{C}$ and $r_2 = 1$. Hence, the pair ($\bm{A}, \bm{B}$) is stabilizable, the pair ($\bm{A}, \bm{C}$) is detectable, and condition (\ref{eq:conditionOfWeightMatrices}) is satisfied, while Newton's and the modified Newton's methods of \cite{ilka2022novel} yield $\bm{K}$ = -0.837. The
two methods returned slightly different $\bm{P}$ matrices, with the
difference between corresponding entries of the two returned
$\bm{P}$ matrices being smaller than 10$^{-7}$, and hence, both methods
gave approximately the same $\bm{G}$ = [0, 0.284, 0.357, 0.089]. Furthermore, condition (\ref{eq:AREoutputFdbk}) was satisfied by both methods
since the left-hand side of the equality contained entries
smaller than 10$^{-7}$. The implementation of the closed-loop system with the initial and set-point conditions used in Section \ref{sub:FullState-CO2case} gives the results reported in Fig. \ref{fig:outputFdbk_CO2}. As can be seen, \emph{the concentration of $\text{CO}_2$ in the target area reaches the value of the pre-industrial era (i.e., 637.2 tonnes) in approximately 60 days}, which is about 35 days slower than the system with full-state feedback control.       
The fact that the action of the output-feedback controller is slower than the full-state one is explained by comparing the control inputs: \emph{$u(0) \approx$ 260 t/d in the former and $u(0) \approx$ 1200 t/d in the latter}, which is motivated by the fact that the full-state feedback uses the knowledge of all four states, whereas output-feedback only one measured state (i.e., $x_1(t)$).  

\begin{figure*}
\centering
\subfloat{
  \includegraphics[width=0.32\textwidth]{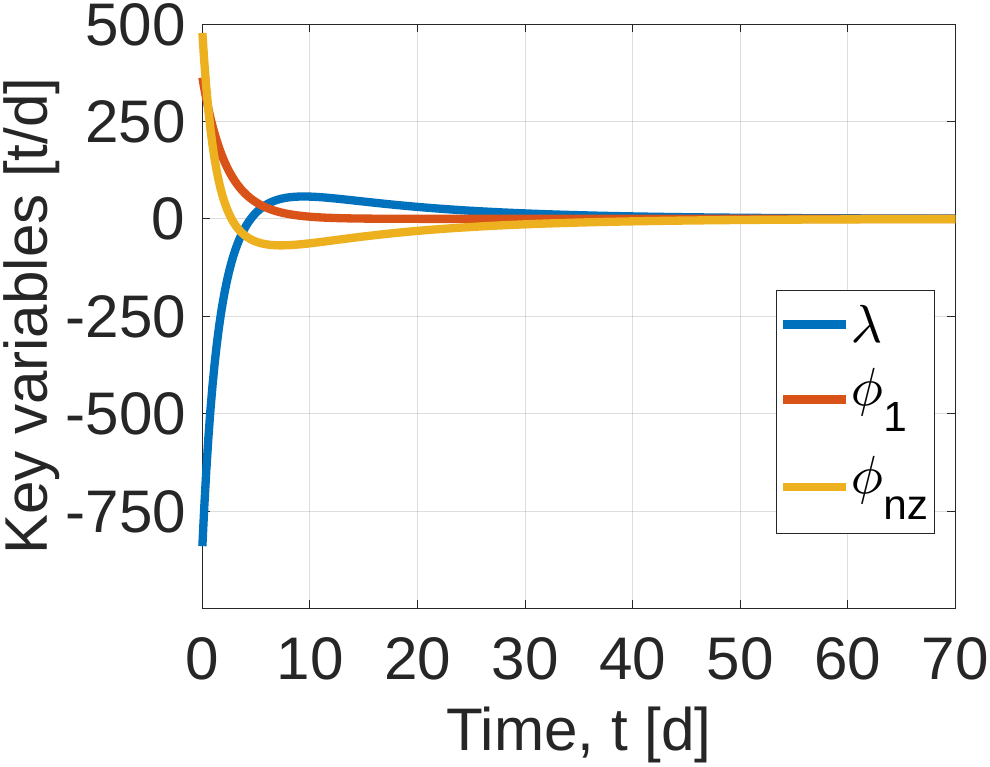}
  }
\subfloat{
  \includegraphics[width=0.27\textwidth]{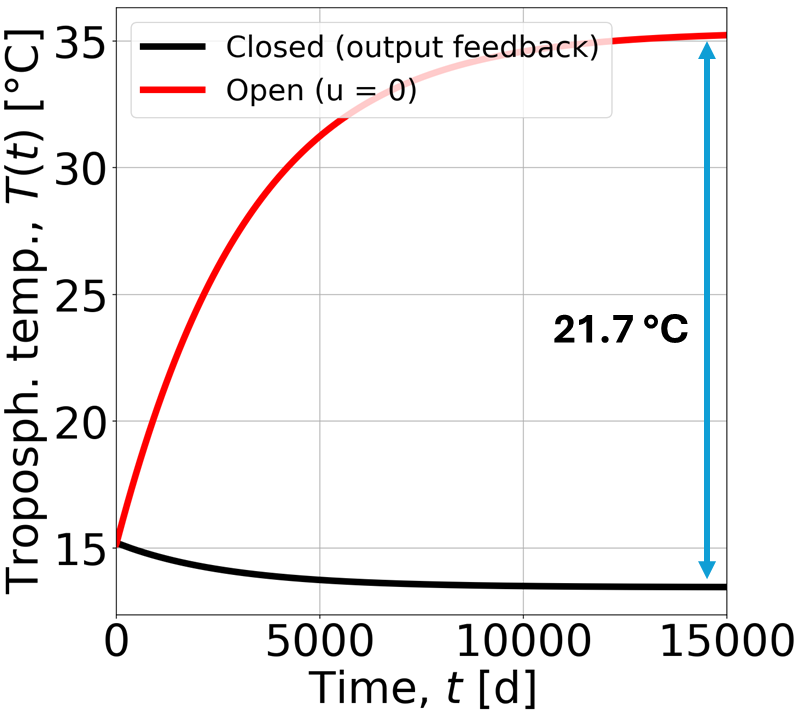}
  }
  \subfloat{
  \includegraphics[width=0.27\textwidth]{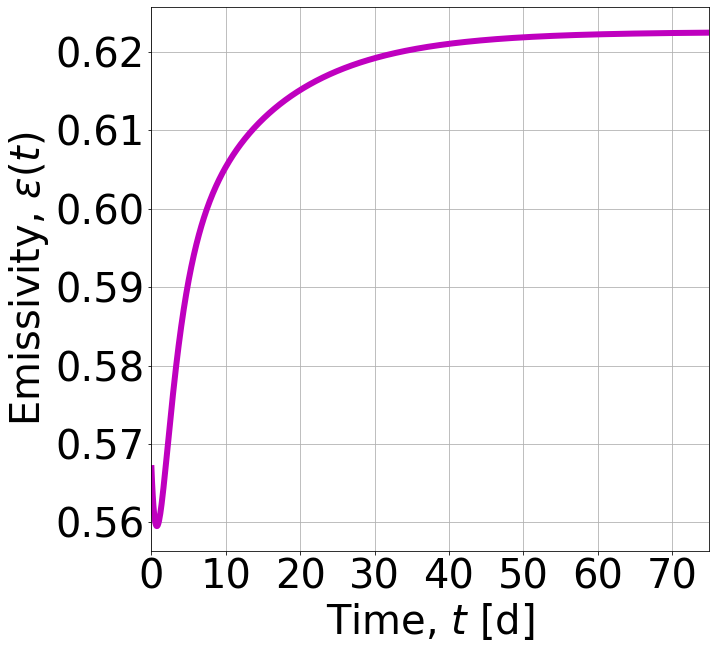}
  }
\caption{Key variables with \emph{output feedback}: (left) circularity, $\phi_1$, and $\phi_\text{nz}$; (middle) tropospheric temperature with open and closed loops; (right) emissivity.}
\label{fig:outputFdbk_CO2_keyVariables}
\end{figure*}

The key variables with output-feedback are shown in Fig. \ref{fig:outputFdbk_CO2_keyVariables}. In particular, $\lambda(0)$ is negative in contrast with the case of full-state feedback since, as discussed above, the control input $u(0)$ is more than 4 times smaller. As in the case of full-state feedback, $\lambda(t)$, $\phi_{1}(t)$, and $\phi_{\text{nz}}(t)$ converge to zero, while the tropospheric temperature behaves analogously to the case of full-state feedback discussed in Section \ref{sub:FullState-CO2case}. In contrast, the emissivity $\epsilon(t)$ has a different response in the full-state case compared to the output-feedback case.    

\section{Conclusion}\label{sec:Conclusion}
This paper proposes CarboNet, a network of thermodynamic compartments controlled via LQRs for tropospheric carbon regulation. The controllers effectively reduced the $\text{CO}_2$ concentration to the pre-industrial era level and also met the net-zero target, but tolerating $\text{CO}_2$ emissions for only 6 days. With the closed-loop control, the tropospheric temperature stabilizes approximately to the pre-industrial era condition, i.e., to 13.5 °C, which is 21.7 °C lower than the steady-state temperature achieved without carbon capture. Future research will focus on controller designs that are able to tolerate the combustion for more than 6 days as well as consider more accurate models of the dynamics of the tropospheric temperature to produce more accurate predictions.




\bibliographystyle{IEEEtran}
\bibliography{References}



\end{document}